\documentclass[a4paper,12pt]{article}

\usepackage{amsfonts}
\usepackage{amsthm}
\usepackage{amsmath}
\usepackage{amsfonts}
\usepackage{latexsym}
\usepackage{amssymb}

\newtheorem{fed}{Definition}[section]
\newtheorem{teo}[fed]{Theorem}
\newtheorem{lem}[fed]{Lemma}
\newtheorem{cor}[fed]{Corollary}
\newtheorem{pro}[fed]{Proposition}
\theoremstyle{definition}
\newtheorem{rem}[fed]{Remark}
\newtheorem{rems}[fed]{Remarks}
\newtheorem{exa}[fed]{Example}
\newtheorem{exas}[fed]{Examples}

\def\bdem{{\bdem\ }\rm}

\newfont{\bb}{msbm10}
\def\Bbb#1{\mbox{\bb #1}}

\def\C{\mathbb {C}}
\def\R{\mathbb {R}}
\def\N{\mathbb {N}}

\def\inc{\subseteq}

\def\EOE{\hfill $\blacktriangle$}
\def\inv{^{-1}}
\def\*A{\#\sb A}

\def\eps{\varepsilon}

\def\H{{\cal H}}

\def\sii{if and only if }

\def\ca{L(\H ) }
\def\cam{L(\H )^+ }

\def\cH{{\cal H}}

\def\cS{{\cal S}}
\def\cT{{\cal T}}

\def\csha{\Sigma (\mathcal{S}, A)}
\def\Csha{\rho (\mathcal{S}, A)}
\def\cshat{\Sigma (\mathcal{S}, A^t)^{1/t}}

\def\rai{^{1/2}}
\def\mrai{^{-1/2}}

\def\api{\langle}
\def\cpi{\rangle}

\def\noi{\noindent}

\def\bm{\left(\begin{array}}
\def\em{\end{array}\right)}
\def\bdem{\begin{proof}}
\def\edem{\end{proof}}
\def\ben{\begin{enumerate}}
\def\een{\end{enumerate}}
\def\beq{\begin{equation}}
\def\eeq{\end{equation}}
\def\barr{\begin{array}}
\def\earr{\end{array}}
\def\inv{^{-1}}

\def\H{{\cal H}}
\def\lh{{L(\H)}}
\def\lh+{{\lh^+}}
\def\cas{{\ca_{sa}}}

\def\la{\lambda}
\def\eps{\varepsilon}

\date{}

\DeclareMathOperator{\leqr}{\preccurlyeq}

\DeclareMathOperator*{\convsotdpre}{\searrow}
\DeclareMathOperator*{\convsotipre}{\nearrow}

\newcommand{\pint}[1]{\displaystyle \left \langle #1 \right\rangle}
\newcommand{\hil}{\mathcal{H}}
\newcommand{\ele}{\mathcal{L}}

\newcommand{\ese}{\mathcal{S}}
\newcommand{\eme}{\mathcal{M}}

\newcommand{\ete}{\mathcal{T}}

\newcommand{\op}{L(\mathcal{H})}
\newcommand{\opsa}{L(\mathcal{H})_{sa}}

\newcommand{\posop}{L(\mathcal{H})^+}

\newcommand{\conv}{\xrightarrow[n\rightarrow\infty]{}}
\newcommand{\convsot}{\xrightarrow[n\rightarrow\infty]{\mbox{\tiny{S.O.T.}}}}
\newcommand{\convsotr}{\xrightarrow[]{\mbox{\tiny{S.O.T.}}}}

\newcommand{\sig}[2]{\rho\left(#1,#2\right)}

\newcommand{\Sig}[2]{\Sigma\left(#1,#2\right)}

\newcommand{\convsotd}{\convsotdpre_{n\rightarrow\infty}^{\mbox{\tiny
\textbf SOT}}}
\newcommand{\convsoti}{\convsotipre_{n\rightarrow\infty}^{\mbox{\tiny
\textbf SOT}}}
\newcommand{\spec}[1]{\sigma\left(#1\right)}
\newcommand{\specr}[1]{\sigma_\rho\left(#1\right)}
\newcommand{\specm}[1]{\sigma_+\left(#1\right)}
\newcommand{\specM}[1]{\sigma_{_-}\left(#1\right)}

\newcommand{\kol}[2]{K\left(#1,#2\right)}
\newcommand{\kkol}[2]{k\left(#1,#2\right)}

\begin{document}


\title{{\Huge\textbf{Spectral shorted operators.}}}
\author{Jorge Antezana, Gustavo Corach and Demetrio Stojanoff}
\maketitle

\centerline{To the memory of Gert K. Pedersen}
\vglue1truecm

\noindent{\bf Jorge Antezana},

\noi Depto. de Matem\'atica, FCE-UNLP,  La Plata,  Argentina 
and IAM-CONICET

\noindent{e-mail: atenzana@mate.unlp.edu.ar}

\noindent{\bf Gustavo Corach},

\noi Depto. de Matem\'atica, Facultad de Ingenier\'\i a UBA,
Buenos Aires, Argentina
and IAM-CONICET

\noindent{e-mail: gcorach@fi.uba.ar}

\medskip

\noindent{\bf Demetrio Stojanoff},

\noi Depto. de Matem\'atica, FCE-UNLP,  La Plata,  Argentina
and IAM-CONICET

\noindent{e-mail: demetrio@mate.unlp.edu.ar}
\medskip

\bigskip

\begin{abstract}{
If   $\mathcal H$ is a Hilbert space, $\mathcal S \subseteq \mathcal H$ 
is a closed subspace of $\mathcal H$, and $A $ is a positive bounded 
linear operator on $\mathcal H$,  the spectral shorted operator
$\rho(\mathcal S,  A)$ is defined as the infimum of the sequence
$\Sigma (\mathcal S, A^n)^{1/n}$, where $\Sigma (\mathcal S, B)$ 
denotes the shorted operator of $B$ to $\mathcal S$. We
characterize the left spectral resolution of $\rho(\mathcal S,  A)$
and show several properties of this operator, particularly in the 
case that $\dim \mathcal S = 1$.
We use these results to generalize the concept of Kolmogorov complexity
for the infinite dimesional case and for non invertible operators.}
\end{abstract}

\bigskip

\noi 

\noi
\date{}

\section{Introduction}

Let $\H$ be a separable Hilbert space and $\ca$ the algebra of
bounded operators on $\H$. Given a positive (i.e. semidefinite 
non negative) operator $A\in \ca$ and a closed subspace $\ese$ 
of $\H$, the shorted operator $\Sig{\ese}{A}$ was defined by Krein \cite{[K]}
and Anderson-Trapp \cite{[AT2]} by 
$$ \csha =
\max \{ X \in \posop : X \le A \quad \hbox { and } \quad R(X)\inc
\cS \} ,
$$ 
where the maximum is taken for the natural order relation in $\posop$, 
the set of positive operators in $\ca$ (see \cite{[AT2]}, \cite{[PEKA]}, \cite{[PS]}).

In a previous paper \cite{[ACS]}, the authors have defined, under the asuumption that 
$\dim \H < \infty$, the so called spectral shorted matrix:   
\beq\label{ladef}
\sig{\ese}{A} = \lim_{m \to \infty} \Sig{\ese}{A^m}^{1/m}
= \inf_{m \to \infty} \Sig{\ese}{A^m}^{1/m}.
\end{equation}
This paper is the continuation of \cite{[ACS]}. It is devoted to study the 
natural generalization of $\rho$ to the infinite dimensional
setting. If $\dim \H = \infty$ and $A \in \cam$, the operator $\sig{\ese}{A}$
is also defined by equation (\ref{ladef}), under the assumption that
the subspace $\ese $ is closed. We call this operator the $spectral \ shorted \ 
operator$ associated to $\ese$ and $A$.

Many properties of the spectral shorted matrices showed in \cite{[ACS]} 
hold also for spectral shorted operators, but some of them must be 
formulated in terms of the spectral measure of $A$ instead of 
eigenvalues and eigenspaces, as in \cite{[ACS]}. 

As in the matrix case, the properties of $\rho$ are strongly related with the so
called $spectral$ $ order$ of positive operators.
Recall the definition of the spectral order $\leqr$ in $\posop$:
given $A,B\in\posop$,
we write  $A\leqr B$ if $A^m\leq B^m$ for all $m\geq 1$.
The spectral order was extensively 
studied by M. P. Olson in \cite{[Ols]}, where the following
characterization is proved: given $A, B \in \posop$, then
$A \leqr B$ \sii  $f(A) \le f(B)$ for every  non-decreasing map
$f : [0 , + \infty ) \to \mathbb R$.

Section 2 contains preliminaries and a brief account of the main properties
of the shorting operation, spectral order and spectral resolutions. 
In section 3 we collect those properties of $\rho$ which can be plainly 
generalized to the infinite dimensional setting.  The 
most subtle tool is the continuity of the map 
$t \mapsto t^r$ (for $0\le r \le 1$) with respect 
to the strong operator topology on $\cam$. 
It is used, for instance, for proving that for every $t >0$,  
\beq\label{alat}
\sig{\ese}{A^t}=\sig{\ese}{A}^t .
\end{equation}
The spectral order provides the following link with
Krein and Anderson-Trapp definition of the shorted operator: $\sig{\ese}{A}$ is the
biggest (in both orders $\le$ and $\leqr$)
element $D$ of $\posop$ such that $D\leqr A$ and
$R(D)\subseteq \ese$ (see  Theorem \ref{caraterizacion uno}).
This shows the monoticity of  $\sig{\ese}{\cdot}$ with respect to the preorder 
$\leqr$ and allows us to get  some  results about limits of spectral shorted
operators. 

In section 3 we get a complete characterization of $\sig{\ese}{A}$ in terms of the 
(left) spectral resolution of $A$: for every $0< \la \in \R$,
$$
\aleph_{[\la , \infty )}  (\sig{\ese}{A}) = \aleph_{[\la , \infty )}  (A) \wedge P_\ese .
$$
This results allows us to get simple proofs in our context of several properties 
of spectral shorted matrices. For example, given $A\in\posop$ and two 
closed subspaces $\ese$ and $\ete$ of $\hil$, 
\ben
\item $\sig{\ese\cap\ete}{A}= \sig{\ete}{\sig{\ese}{A}} .$
\item $\spec{\sig{\ese}{A}}\subseteq \spec{A}$.
\item $f(\sig{\ese}{A})=\sig{\ese}{f(A)}$, for every non-decreasing $right$
$continuous$ positive function $f$ defined on $[0,+\infty)$.
\item $\la_{min} (A) P_\ese \le \sig{\ese}{A}$, where $\la_{min} (C) = \min \spec{C}$, for 
$C \in \cam$.
\item If $\sig{\ese}{A}$ is considered as acting in $\ese$, then 
$$\la_{min}(\sig{\ese}{A})
 = \min \{\mu\in\spec{A}:\; P_\ese \  \aleph _{[\mu , \mu
+ \eps ) } (A)  \not = 0\;\;\forall \; \varepsilon>0\}.
$$
\een
The case $\dim \ese = 1$ is extensively studied in section 5. 
If $\ese$ is the subspace generalted by 
the unit vector $\xi$, we denote by $\sig{A}{\xi}$ 
the unique positive number such that 
$
\sig{\ese}{A} = \sig{A}{\xi} P_\ese .
$
The main results of this secton are: 
\ben
\item If $A \in \cam$ and $\xi\in \H $ is an unit vector, then 
$$
\sig{A}{\xi} = \min \spec{\sig{\ese}{A}} 
 = \min \Big\{ \mu\in\spec{A}: \;
\aleph _{[\mu , \mu + \eps ) } (A)\xi  \not = 0 \; \; \forall \; \varepsilon>0 \Big\}
.
$$
\item $\sig{A}{\xi} = \max \{ \la \in \spec{A} : \xi \in R(\aleph_{[\la , \infty )} (A)) \}.$
\item If $A$ is invertible, then $\displaystyle 
\sig{A}{\xi} = \lim _{m \to \infty } \|A^{-m} \xi \|
^{-1/m} = \inf _{m \in \Bbb N}  \|A^{-m} \xi \| ^{-1/m} .
$
\item 
If $R(A)$ is closed and $\xi \in R(A)$, then, 
$\sig{A}{\xi} = \lim _{m \to \infty } \|(A^\dagger)^m \xi \| ^{-1/m}$, 
where $A^\dagger$ is the Moore-Penrose 
pseudo-inverse of $A$. If $\xi \notin R(A)$, then  $\sig{A}{\xi} =0$.
\item
If $\specr{A} = \big\{ \sig{A}{\xi} : \|\xi \| = 1\big\}$, then 
$$\specr{A} = \specm{A} \cup \sigma_{pt}(A)
=\big\{ \la \in
\spec{A}: \forall \ \eps >0 \ , \ \aleph_{[\la, \la+ \eps )}(A) \neq 0
\},
$$
where $\sigma_{pt}(A)$ denotes the point spectrum of $A$, i.e the set of eigenvalues of $A$
and $\specm{A}$ is the set of points in $\spec{A}$ which are limit point of
$\spec{A} \setminus \{\la\}$ from the right. This shows that $\specr{A}$ is allways
dense in $\spec{A}$, but $\specr{A} \neq \spec{A} $ in general.
\item $\sig{A}{\xi} \neq 0$ \sii $\xi \in R_0(A):=
\bigcup_{\la > 0} R(\aleph_{[\la , \infty)} (A)) \inc R(A). $
\een

\noi
In \cite{[FF]}, J. I. Fujii and M. Fujii consider the Kolmogorov's
complexity
\beq\label{FF}
\kol{A}{\xi} = \lim _{n \to \infty} \frac{\log (\pint {A^n \xi , \xi})}{n} 
=\log \lim _{n \to \infty} \pint {A^n \xi , \xi}  ^{1/n} .
\end{equation}
for an invertible positive matrix $A$ and a  unit vector $\xi$ and
show several properties $K$. In \cite{[ACS]} we show that, if $\ese$ is the
subspace generated by $\xi$, then
$$
\kol{A}{\xi} = \log \sig{A^{-1}}{\xi} \inv .
$$
We define a generalized version (for  $\dim \H = \infty$ and $A\in \cam$ not 
necessarily invertible) 
without logarithms (in order to avoid the value $-\infty$) 
of the Kolmogorov complexity  as follows:
given $\xi \in \H$, $\xi \neq 0$  and $A\in \cam$, we denote by
$$
\kkol{A}{\xi} = \lim _{n \to \infty} \pint {A^n \xi , \xi}  ^{1/n} ,
$$
that is, $\kkol{A}{\xi} = \exp \kol{A}{\xi}$ in the cases where
$\kol{A}{\xi}$ can be defined as in equation (\ref{FF}).
Among other properties, we show that: if $\xi \in \H$ and $A\in \cam$, then
\ben
\item If $\|\xi\| = 1$, then the sequence $\pint {A^n \xi , \xi}  ^{1/n}$ is increasing.
So that, $\lim _{n \to \infty} \pint {A^n \xi , \xi}  ^{1/n}$  exists
for every $\xi \in \H$.
\item $\kkol{A}{\xi} = \kkol{A}{ a \xi}$ for every $0 \neq a \in \C$.
\item $\kkol{A}{\xi} = \kkol{A}{\aleph_{[\la , \infty)} (A) \xi}$ for every
$\la >0 $ such that $\aleph_{[\la , \infty)} (A)\xi \neq 0$.
\item $\kkol{A}{\xi} \neq 0$ (i.e. $\kol{A}{\xi} \neq -\infty$)
\sii $P_{\overline{R(A)}}\  \xi \in R_0(A)\setminus \{ 0\}$.
\item If $\xi \neq 0$, then $\kkol {A}{\xi}\in \spec{A}$. Moreover,  
$$
\big\{  \kkol {A}{\xi} : \xi \neq 0\} =
 \big\{ \la \in \spec{A}:
 \aleph_{(\la+ \eps, \la ]}(A) \neq 0 \ , \ \forall \ \eps >0 \ \} , 
$$  
which is a dense subset of $\spec{A}$. 

\item \  \\ 
$
\barr{rl}
\kkol{A}{\xi} &= \min \big\{ \la \in \spec{A} : \xi \in R(\aleph_{(-\infty , \la]} (A))\big\}
\\&\\
& = \max \Big\{ \mu\in\spec{A}: \;
\aleph _{( \mu - \eps, \mu] } (A)\xi  \not = 0 \; \; \forall \; \varepsilon>0 \Big\} 
\\&\\
& = \sup \Big\{ \mu\in\spec{A}: \; \aleph_{[\mu , \infty)} (A)\xi \neq 0 \Big\} \ .
\earr
$
\item If $R(A)$ is closed, then 
\ben
\item If $\xi \in R(A)$ then $\kkol{A}{\xi} = \sig{A^\dag}{\xi}\inv$.
\item If $\xi \notin  R  ( A)$, but $P\xi \neq 0$, where $P = P_{R(A)}$,  then
$$
\kkol{A}{\xi} = 
 \kkol{A}{P\xi} =  \sig {A^\dag }{\frac{P\xi}{\|P\xi\|}}\inv  \ .
 $$
\een
\een

\section{Preliminaries}

For an operator $A \in \op$, we denote by $R(A)$ the range of $A$,
$ N  ( A)$ the null-space of $A$, $\sigma (A)$ the spectrum of $A$,
$A^*$ the adjoint  of $A$, $\rho(A)$ the spectral radius of
$A$, $\|A\|$ the spectral norm (i.e. the operator norm induced by
the norm of the Hilbert space $\hil$) of $A$.
We denote by $\cas$ the space of selfadjoint
operators in $\ca$ and by $\cam$ the space of positive (i.e. semidefinite
non-negative) operators in $\ca$. If $A \in \cas$, we
denote by $\la_{min}(A) = \min \sigma (A) = \inf_{\|\xi \| = 1}
\api A\xi , \xi \cpi$.

Given a closed subspace $\ese$ of $\hil$, we denote by $P_\ese$
the orthogonal (i.e. selfadjoint) projection onto $\ese$. If $P$ and 
$Q$ are orthogonal projections, we denote by $P\wedge Q$ the orthogonal 
projection  onto $R(P) \cap R(Q)$. If $B
\in \op $ satisfies $P_\ese B P_\ese = B$, we sometimes consider the
compression of $B$ to $\ese$, (i.e. the restriction of $B$ to
$\ese$ as a linear transformation form $\ese$ to $\ese$), and we
say that we consider $B$ as $acting$ on $\ese$. Several times this is
done in order to consider $\sigma (B)$ just in terms of the action
of $B$ on $\ese$. For example, if $B \ge \la P_\ese$ for some $\la
>0$, then we can deduce that $0 \notin \sigma (B)$, if we consider
$B$ as acting on $\ese$.

Along this note we use the fact that every closed subspace $\ese$
of $ \hil$ induces a representation of elements of $\op$ by $2
\times 2$ block matrices, that is, we shall identify each
$A\in\op$ with a $2\times 2$-matrix, let us say
$
\begin{pmatrix}
  A_{11} & A_{12} \\
  A_{21} & A_{22}
\end{pmatrix}\begin{array}{cc}
  \ese  \\
  \ese^\bot
\end{array} .
$ Observe that $
\begin{pmatrix}
  A_{11}^* & A_{21}^* \\
  A_{12}^* & A_{22}^*
\end{pmatrix}$ is the matrix which represents $A^*$.

We use in this note several standard results of spectral theory, functional calculus
and weak convergences of opeartors in $\cas$. About these matters,
we refer the reader to the books of Pedersen \cite{[Ped2]} or Kadison and Ringrose
\cite{[KR]}.  If $A\in \cas$ we denote by $E_A$ the spectral measure associated
to $A$, defined by $E_A(\Delta ) = \aleph_\Delta (A)$, for any Borel set $\Delta
\inc \R$. The sigles SOT are used to mention the $strong \ operator \ topology $
of $\cas$. In the following subsections, we state explicitely several known 
results which we shall need in the sequel. Particulrly those we think are not 
"de libro".

\subsection*{Shorted operators.}

Following Anderson and Trapp \cite{[AT1]}, \cite{[AT2]}, we define
\begin{fed}\label{definicion de Anderson Trapp} \rm
 Let $A\in \posop$ and $\ese $ a closed subspace of $\hil$. Then,
the $shorted$ $operator$ of $A$ to $\ese$ is defined by 
$$ 
\csha = \max \{ X \in \posop : X \le A \quad \hbox { and } \quad R(X)\inc
\cS \}, 
$$ 
where the maximum is taken for the natural order
relation in $\posop$ (see \cite{[AT2]}).
\end{fed}

\noi 
In the next theorem we state some results on shorted operators
proved by Anderson and Trapp \cite{[AT2]},  M.G. Krein \cite{[K]}
and E. L. Pekarev \cite{[PEKA]} which are relevant in this paper.

\begin{teo}\label{propiedades del shorted}
Let $\ese$ and $\ete$ be subspaces of $\hil$ and let $A,B \in
\posop$. Then
\ben
  \item If $\ese\subseteq \ete$, then, $\Sig{\ese}{A}\leq
  \Sig{\ete}{A}$.
  \item $\Sig{\ese\cap\ete}{A}=\Sig{\ese}{\Sig{\ete}{A}}$.
  \item If $A\leq B$, then, $\Sig{\ese}{A}\leq \Sig{\ese}{B}$.
\item Let $\eme = A\mrai (\ese )$. Then  $\Sig{\ese}{A}= A\rai P_\eme A\rai$.
\een
\end{teo}

\noi There are also some results about the continuity of the
shorting operation (see \cite{[AT2]}, Corollary 2 and 3).

\begin{pro}\label{convergencia sot del shorted}
Let $A_n$ ($n\in\mathbb{N}$) be a sequence of positive matrices
such that $\displaystyle A_n\convsotd A$. Then, for every closed
subspace $\ese$ it holds
$
\displaystyle \Sig{\ese}{A_n}\convsotd\Sig{\ese}{A}.
$
\end{pro}

\begin{pro}\label{shur abajo}
Let $\ese_n$ ($n\in\mathbb{N}$) and $\ese$ be closed subspaces
such that $\displaystyle P_{\ese_n}\convsotd  P_{\ese}$. Then, for
every $A\in \posop$, it holds that 
$
\displaystyle \Sig{\ese_n}{A}\convsotd \Sig{\ese}{A}.
$
\end{pro}
\bdem Since $\{\Sig{\ese_n}{A}\}$ is a non-increasing sequence, it
has a strong limit, say $L$. As $\Sig{\ese_n}{A}\leq A$ for all
$n\in\mathbb{N}$, then $L\leq A$. On the other hand, $L\le \Sig{\ese_n }{A}$ 
implies
\[
R(L \rai )\subseteq R\left(\Sig{\ese_n}{A} \rai \right)\subseteq
\ese_n \hspace{0.5cm} \forall\;n\in\mathbb{N} .
\]
Therefore  $R(L)\subset
\bigcap^\infty_{n=1}\ese_n = \ese $. Finally, if $0\leq X\leq A$
and $R(X)\subset \ese$, then $R(X)\subseteq \ese_n$, so that
$X\leq \Sig{\ese_n}{A}$, for all $n\in\mathbb{N}$. Therefore
$X\leq L$. \edem

\subsection*{Spectral order.}


The spectral order was considered by Olson (see \cite{[Ols]}) with
the purpose of reporting an order relation with respect to which
the real vector space of selfadjoint operators form a
conditionally complete lattice. Throughout this note we shall only
use the spectral order for positive operators, and this is the
reason why we take the following statement as definition of the
spectral order.

\begin{fed} \rm \label{spectral order}
Let $A, B \in \posop$. We write $A\leqr B$ if for every
$m\in\mathbb{N}$ it holds that $A^m\leq B^m$. The relation $\leqr$
defined on $\posop$  is a partial order and it is called
$spectral$ $order$.
\end{fed}

\begin{exas} Given $A, B \in \posop$. Then
\ben
  \item If $AB=BA$ and $A\leq B$, then, $A\leqr B$.
  \item If $\dim \cH = n < \infty$, then $A\leqr B$ \sii
  there is a positive integer $k\leq n$ and an sequence of
    positive matrices $\{D_i\}_{0\leq i\leq k}$ such that,
  $D_0=A$, $D_k=B$, $D_i \le D_{i+1}$ and $D_iD_{i+1}=D_{i+1}D_i$ ($i=0,\cdots, k-1$)
  (see \cite{[ACS]}). \EOE
  \een
\end{exas}

\noi
The next results was proved by Olson in \cite{[Ols]}.

\begin{teo}\label{olson}
Let $A,B\in\posop$. The following statements are mutually
equivalent.
\begin{description}
  \item[(1)] $A\leqr B$,
  \item[(2)] $\aleph_{[\la , \,\infty )}(A)\leq
\aleph_{[\la ,\, \infty
  )}(B)$ ($0\leq \la<\infty$),
  \item[(3)] $f(A)\leq f(B)$ for every non-decreasing
continuous function  $f$ on $[0,\infty)$.
\end{description}
\end{teo}

\medskip
\noi
The following result about functions which are continuous relative to the
S.O.T topology of $\cam$ or $\opsa$ is a key tool for the extention of the results
about spectral shorted operators from matrices to operators in Hilbest spaces.
A proof can be found, for example, in Pedersen's book
\cite {[Ped]}.

\medskip
\begin {lem}\label{xalar}
Let $f: \R \to \R$ be a continuous function such that $f(0) = 0$ and
$|f(t)| \le \alpha |t| + \beta$ for some positive numbers
$\alpha $ and $\beta$.
Then, if  $\{A_\alpha\}_{\alpha \in \Lambda} $ is a net in $\opsa$ such that
$\displaystyle A_\alpha\convsotr A \in \opsa$, it holds that
$\displaystyle f(A_\alpha)\convsotr f(A) $, i.e. $f: \opsa \to \opsa$ is continuous for the
S.O.T. topology.
In particular $f(t) = t^r $ for $0\le r\le 1$ is S.O.T.-continuous
in $\cam$.
\end{lem}

\smallskip
\begin{pro} \label {abajo}
Let $\{A_n\}$ be a sequence in $\cam$ such that $A_{n+1}\leqr A_n$,
$n \in \mathbb N$ and $\displaystyle A_n\convsotd A \in \cam$.
Then, for every $k \in \N$, $\displaystyle A_n^k \convsotd A^k $. In particular,
$A \leqr A_n $, $n \in \N$.
\end{pro}
\bdem Fix $k \in \N$. Since the sequence $\{A_n\}$ is non increasing with respect to the
spectral order, there exists  $B \in \cam$ such that $\displaystyle A_n^k \convsotd B$.
By Lemma \ref{xalar}, applied to the map $f(t) = t^{1/k}$,
we can deduce that $\displaystyle A_n\convsotd B ^{1/k} = A$. So that,
$B = A^k$. \edem

\subsection*{Spectral resolutions}
Given $f : \mathbb R \to \ca$, we say that $f$ is a right (resp. left)
$spectral \ resolution$ if
\ben
\item There exist $m, M \in \R$ such that $f(\la ) = 0$
for $ \la <m$ and $f(\la ) = I$ for $\la >M$ (resp. $f(\la ) = I$
for $ \la <m$ and $f(\la ) = 0$ for $\la >M$).
\item $f(\la )$ is a selfadjoint projection, for every $\la \in \R$.
\item If $\la\le \mu$ then $f(\la ) \le f(\mu ) $ (resp. $f(\la ) \ge f(\mu ) $) as operators.
\item $f$ is continuous on the right (resp. $f$ is continuous on the left).
\een
Under these hypothesis, by the standard spectral theory,
there exists an unique $A  \in \cas$ such that $f$ is its spectral
resolution, i.e.
\beq\label{sr}
f(\la ) = E_A ( \ (-\infty , \la ]) = \aleph_{(-\infty , \la ]\ }(A) \quad
\hbox { (resp. } \
f(\la ) = E_A ( \ [\la, \infty)\ ) = \aleph_{[\la, \infty)}(A) \ ).
\end{equation}
Conversely, if $A  \in \cas$, then the map $f$ defined by equation (\ref{sr})
is a right (resp. left) spectral resolution.

The relation between right and left spectral resolutions is given by the following identity: if
$A \in \cas$, then $E_A ( \ [- \la, \infty)\ ) = E_{-A} ( \ (-\infty ,  \la ])$. On the other hand,
if $f$ is a left spectral resolution, then $g(\la ) = f(- \la )$ is a right spectral resolution.
Then, if $A$ is the operator associated to $g$, then $-A$ is the operator associated to $f$.

\section{The spectral shorted operator}


In this section we define the spectral shorted operator in the infinite 
dimensional setting, and we show its basics properties. All results and proofs 
of this section ar very similar as those which appear in \cite {[ACS]} for 
de finite dimensional case, but using SOT-convergence instead of convergence in norm. 
The main difference is that, in the proof of 
Proposition \ref{the key}, we need to apply Lemma \ref{xalar} about SOT-continuity
of the map $A \mapsto A^r$ for $0\le r \le 1$. Also Propostion \ref{convegencia sot} 
is a properly infinite dimensional result.

\medskip
\begin{pro}\label{schat} Let $A \in \ca^+ $ and
$\cS\subseteq\hil$ a closed subspace. Then the map 
$t \mapsto  \cshat , \ t\in [1, \infty ) $  is non-increasing.
\end{pro}
\bdem Fix $t \ge 1$. Then $\Sigma (\mathcal{S}, A^t) \le A^t$. Since $0 \le
1/t \le 1$, by L\"owner theorem we can deduce that $\cshat \le A .$ On the other hand
$R(\cshat ) \inc \cS$. So, by the
definition of shorted operator, $ \cshat \le \csha$. Now, given $1\le r \le s $, take 
$t = s/r\ge 1$. By the previous remarks, applied to $A^r$ and $t$, we have that 
$$
\Sig{\ese}{ A^r } \ge  \Sig{\ese}{ A^{rt}}^{1/t} = \Sig{\ese}{ A^s}^{r/s} .
$$
Since $1/r \le 1$, by L\"owner theorem we have that
$\Sig{\ese}{ A^r } ^{1/r}\ge  \Sig{\ese}{ A^s}^{1/s} $.
\edem

\begin{fed} \rm Let $A \in \ca^+ $, 
and let $\ese \inc \H$ be a closed subspace.
We denote by
$$
\Csha = \inf_{t\ge 1} \cshat  = \lim_{t \to + \infty } \  \cshat  ,
$$
where the limit is taken in the strong operator topology (SOT).
\end{fed}

\begin{rem}\rm\label{ejem}
Let $A\in\posop$ and let $\ese$ and $\ete$ be closed subspaces. 
\begin{enumerate}
  \item If $A =P_\ete$,   then 
$\Csha = \cshat =  P_{\cS \cap \cT} $, for every $t\in [1,\infty)$.
  \item If $AP = PA $, then $\Csha = \cshat = P A $, 
for every $t\in [1, \infty)$.
  \item $\sig{\ese}{cA}=c\,\sig{\ese}{A}$ for every $c\in [0,+\infty)$.
  \item If $\ese\subseteq \ete$, then, $\sig{\ese}{A}\leq\sig{\ete}{A}$, since 
 $\Sig{\ese}{A^t}^{1/t} \le \Sig{\ete}{A^t}^{1/t}$ for every $t\ge 1$. 
  \EOE
\end{enumerate}
\end{rem}


\begin{pro}\label{the key}
Let $A\in\posop$ and $\ese\inc\hil$ be a closed subspace. Then, for
every $t\in(0,\infty)$ it holds that
\[
\sig{\ese}{A}^t=\sig{\ese}{A^t}
\]
In particular, for every $t\in(0,\infty)$
\[
\sig{\ese}{A}^t\leq A^t
\]
\end{pro}
\bdem
Firstly, we prove the statement for $t\geq 1$.
By Lemma \ref{xalar}, the map $x\rightarrow x^r$ is continuous in the strong
operator topology  when $0\leq r\leq 1$.
So, given $t\in(1,\infty)$, since $st\rightarrow \infty$ as $s\rightarrow
\infty$, we have that
\begin{align*}
\sig{\ese}{A^t}^{1/t}&=\left(\lim_{s\rightarrow \infty}
\Sig{\ese}{(A^t)^s}^{1/s} \right)^{1/t}=\lim_{s\rightarrow
\infty} \Sig{\ese}{A^{st}}^{1/st}=\sig{\ese}{A},
\end{align*}
where the limits are taken in the strong operator topology.
This proves, for $t\ge1$, that
\begin{equation}\label{potencias para t mayor que uno}
\sig{\ese}{A^t}=\sig{\ese}{A}^t.
\end{equation}
Now, if $t\in(0,1)$,
\[
\sig{\ese}{A^t}=\left(\sig{\ese}{A^t}^{1/t}\right)^t=
\sig{\ese}{(A^t)^{1/t}}^t=\sig{\ese}{A}^t,
\]
where in the second equality, we have used equation (\ref{potencias para t
mayor que uno}) for $\displaystyle \frac{1}{t}\geq1$.\edem


\bigskip
\noi
Recall that given two positive operators $A$ and $B$ we say that
\begin{align*}
A\leqr B \hspace{1cm}\mbox { if } \hspace{1cm} A^n\leq B^n\;\;\;\forall n\geq
1
\end{align*}
With respect to this order, the spectral shorted operator has a
characterization similar to Anderson-Trapp's definition of shorted
operator.

\begin{teo}\label{caraterizacion uno}
Let $A\in\posop$ and $\ese$ a closed subspace of $\hil$. If
\begin{align*}
\eme_\rho(\ese,A)= \ \{D\in\posop:\;\;D\leqr A,\;\; R(D)\subseteq
\ese\}
\end{align*}
then
\begin{align*}
\sig{\ese}{A}=\max\eme_\rho(\ese,A),
\end{align*}
where the ``maximum" is taken  for  any of the orders $\le$ and
$\leqr$.
\end{teo}
\bdem Firstly, note that $\sig{\ese}{A}\in\eme_\rho(\ese,A)$. In
fact, $\sig{\ese}{A}^m\leq A^m$ for every $m\in \mathbb{N}$ by
Proposition \ref{the key}, and 
$R(\sig{\ese}{A})\subseteq\ese$ by definition.

Suppose that $D\in\eme_\rho(\ese,A)$. Fix $m \in \N$. As $D^m\leq A^m$, it
holds that 
$
\Sig{\ese}{D^m}^{1/m}\leq \Sig{\ese}{A^m}^{1/m}.
$
Since $\Sig{\ese}{D^m}^{1/m}=D$, 
taking $m \to \infty$ we have
$D\leq \sig{\ese}{A}.$ This shows that 
$\sig{\ese}{A}=\max\eme_\rho(\ese,A)$ for the usual order.

Note also that, if $D \in \eme_\rho(\ese,A)$, then for every $k
\in \mathbb N$, $D^k \leqr A^k$ and $D^k \in \eme_\rho(\ese,A^k)$. 
By the previous case, applied to $A^k$, one gets
\[
 D^k \le \sig{\ese}{A^k} = \sig{\ese}{A}^k \ , \quad k\in \N .
\]
Hence $D\leqr \sig{\ese}{A}$.\edem

\begin{cor}\label{monotonia}
Let $A$ and $B$ be positive operators such that $A\leqr B$ and
$\ese$ and $\ete$ be closed subspaces such that $\ese\subseteq
\ete$. Then $\sig{\ese}{A}\leqr \sig{\ete}{B}$.
\end{cor}
\bdem It is enough to note that $\eme_\sigma(\ese,A)\subseteq
\eme_\sigma(\ete,B)$.\edem


\noi
Another application of  Theorem \ref{caraterizacion uno} is the
following result about the convergence of sequences of spectral
shorted operators.

\begin{pro}\label{convegencia sot}
Let $\{A_n\}$ be a sequence in $\cam$ such that $A_{n+1}\leqr
A_n$, $n \in \mathbb N$ and $\displaystyle A_n\convsot A$, and let
$\{\ese_n\}$ be a sequence of subspaces such that
$\ese_{n+1}\subseteq \ese_{n}$. Then
\[
\sig{\ese_n}{A_n}\convsotd\sig{\ese}{A},
\]
where $\displaystyle \ese=\bigcap_{n=1}^\infty \ese_n$.
\end{pro}
\bdem
By Corollary \ref{monotonia}, for every $n \in \N$, 
$\sig{\ese_{n+1}}{A_{n+1}}\leq \sig{\ese_n}{A_n}$. 
Then there is a positive operator $L$ such
that $\sig{\ese_n}{A_n}\convsot L$. On one hand, by Proposition
\ref{abajo}, $A \leqr A_n$, $n \in \N$. 
As, in addition, $\ese\subseteq \ese_n$,
we have that $\sig{\ese}{A}\leq \sig{\ese_n}{A_n} $, $n \in \mathbb N$. 
This shows that $\sig{\ese}{A}\leq L $.
On the other hand, for every $n>m$ and $k\geq 1$, by Corollary
\ref{monotonia} and the definition of spectral shorted operators, 
\begin{equation}\label{dos del limite comun}
L\leq \sig{\ese_n}{A_n}\leq \sig{\ese_m}{A_n}\leq
\Sig{\ese_m}{A_n^k}^{1/k} .
\end{equation}
Now fix $k\geq 1$. By Proposition \ref{abajo}, $\displaystyle
A_n^k\convsotd A^k$. Therefore, by Lemma \ref{xalar}, 
\begin{equation}\label{tres del limite comun}
\Sig{\ese_m}{A_n^k}^{1/k}\convsotd \Sig{\ese_m}{A^k}^{1/k} .
\end{equation}
In a similar way, using Proposition \ref{shur abajo}, we have that
\begin{equation}\label{cuatro del limite comun}
\Sig{\ese_n}{A^k}^{1/k}\convsotd \Sig{\ese}{A^k}^{1/k}.
\end{equation}
Hence, joining equations (\ref{dos del limite comun}) (\ref{tres
del limite comun}) and (\ref{cuatro del limite comun}), we obtain
$
L\leq \Sig{\ese}{A^k}^{1/k}.
$
Finally, since the last inequality is true for every $k$, by
taking limit we have that $L\leq \sig{\ese}{A}$.\edem

\bigskip
As the following example shows, the last Proposition fails, in
general, if the sequence of subspaces is not non-increasing.

\begin{exa}
Let $\hil$ be a separable Hilbert space and $A$ a positive and
injective operator such that $R(A^{1/2})\neq\hil$. Let $\ele$ be a
proper dense subspace of $\hil$ such that
$R(A^{1/2})\cap\ele=\{0\}$. Take an orthonormal basis $\{e_n\}$ of $\hil$
contained in $\ele$, and define $\ese_n=\langle
e_1,\ldots,e_n\rangle$. Then, $\displaystyle P_{\ese_n}\convsoti
I$, but, $\sig{\ese_n}{A}=\Sig{\ese_n}{A}=0$ for all
$n\in\mathbb{N}$, because, as it was proved in \cite{[AT2]}, 
$R(\Sig{\ese_n}{A}\rai ) = R(A\rai) \cap \ese_n = \{0\}$.
\EOE
\end{exa}

\section{Main properties of $\sig{\ese}{A}$.}

Let $A\in\posop$ and let $\ese$ be a closed subspace of $\H$. It is shown in \cite{[ACS]}
that, if $\dim \cH < \infty$ and $0< \la \in \R$, then
$$ 
\bigoplus _{\mu \ge \la } \ker ( \Csha - \mu I) = \ese \cap
\bigoplus _{\mu \ge \la } \ker ( A - \mu I) .
$$
This can be reformulated, in terms of spactral measures, as 
$$
\aleph_{[\la , \infty )}  (\sig{\ese}{A}) = \aleph_{[\la , \infty )}  (A) \wedge P_\ese .
$$
This formula, which allows to compute the spectrum and  the eigenvectors of $\sig{\ese}{A}$, gives
the complete characterization of $\sig{\ese}{A}$ in the matrix case.

In the infinite dimensional case, a similar formula can be proved following the same methods
(with considerable more effort). Instead, it seems more convenient to construct an
operator by means of  the left spectral resolution given by
\beq\label{SR}
f(\la ) = \left\{ \barr{cl}   \aleph_{[\la , \infty )}  (A) \wedge P_ \ese   & \la > 0 \\
                                 I                  & \la \le 0      \earr \right.
\end{equation}
and then to show that its associated operator agrees with $\sig{\ese}{A}$.
This can be done by using the characterization of $\sig{\ese}{A}$ given in
 Theorem \ref{caraterizacion uno}. Note that the verification of the fact that 
 $f$ is, indeed, a left spectral resolution is apparent from the fact that $\la \mapsto 
 \aleph_{[\la , \infty )}  (A)$ is the left spectral resolution of $A$.

\begin{teo}\label{trunquetti}
\rm Let $A\in\posop$ and let $\ese$ be a closed subspace of $\H$.
Then $\sig{\ese}{A}$ is the operator defined by the
left spectral resolution $f$ defined in equation (\ref{SR}).
In other words, for $0< \la \in \R$,
$$
\aleph_{[\la , \infty )}  (\sig{\ese}{A}) = \aleph_{[\la , \infty )}  (A) \wedge P_\ese .
$$
\end{teo}
\bdem
Let $B$ be  the operator
defined by the spectral resolution $f$. By Theorem \ref{olson}, it is clear that
$B \leqr A$ and every $D \in \eme_\rho(\ese,A)$ satisfies $D\leqr B$. Indeed, suppose that
$0\le D \leqr A $ and $R(D) \inc \ese$. Then, for $\la > 0$,  $\aleph_{[\la , \infty )} (D)
\le \aleph_{[\la , \infty )} (A) $ and
$$
\aleph_{[\la , \infty )} (D) \le \aleph_{(0 , \infty )}(D) \le P_{\overline{R(D)}} \le   P_\ese.
$$
Therefore $\aleph_{[\la , \infty )} (D) \le  \aleph_{[\la , \infty )}  (A) \wedge P_ \ese =
\aleph_{[\la , \infty )} (B)$. Since $\aleph_{[\la , \infty )} (D) = I = \aleph_{[\la , \infty )} (B)$
for $\la \le 0$, we get that $D\leqr B$ by Theorem \ref{olson}.
Finally, since
$$
\aleph_{[\la , \infty )}( \| A\| P_\ese) = \left\{ \barr{ccl}   0 && \|A\| <\la  \\
                                                     P_ \ese   && 0< \la \le \|A\|   \\
                                 I                  && \la \le 0      \earr \right. ,
$$
we deduce that $B \leqr \|A\| \ P_S$ and, in particular, $R(B) \inc \ese$.
Then, by  Theorem \ref{caraterizacion uno},
$$
B = \max \eme_\rho(\ese,A)= \sig{\ese}{A}.
$$
\edem

\begin{cor}\label{syt}
Let $A\in\posop$ and let $\ese$ and $\ete$ be closed subspaces of $\H$.
Then
\begin{align*}
\sig{\ese\cap\ete}{A}= \sig{\ete}{\sig{\ese}{A}} .
\end{align*}
\end{cor}
\bdem
Note that both operators have, as left spectral resolution, the map
$$
f(\la ) = \left\{ \barr{cl}   
\aleph_{[\la , \infty )}  (A) \wedge P_ \ese \wedge P_ \ete  & \la > 0 \\
                                 I                  & \la \le 0      \earr \right.  .
$$
\edem

\begin{rem}\rm
Let $A\in\posop$ and let $\ese$ and $\ete$ be closed subspaces of $\H$.
Then $$\sig{\ese\cap\ete}{A}\leq \sig{\ete}{\Sig{\ese}{A}}.
$$
Indeed, 
it can be deduced from  inequalities
\begin{align*}
\Sig{\ese\cap\ete}{A^{2^m}}\leq
\Sig{\ete}{\Sig{\ese}{A^{2^m}}}\leq
\Sig{\ete}{\Sig{\ese}{A}^{2^m}} &&\forall m\in\mathbb{N} . 
\end{align*}
Note that the mentioned statement can not be deduced from Corollary \ref{syt}.
\EOE

\end{rem}

\bigskip

\medskip
\begin{pro}\label{cotita}
If $A \in \posop$ and  $\mu = \min \spec{A}$, then
$$
\mu P \le \Csha . $$ In particular, if $A$ is invertible  then
$\Csha $ is invertible if it is considered as
acting on $\ese$.
\end{pro}
\bdem Note that $\mu ^ m = \min \spec{A^m}$ for all $m \in \Bbb N$.
Then  $\mu^m  P_\ese \le  \mu^m  I \le A^m $ for all $m \in \Bbb N$. So
that, $\mu P_\ese \leqr A $ and the result follows by  Theorem \ref{caraterizacion uno}. 
\edem

\bigskip
\begin{rem} \label {keeping}\rm
Given an operator $A\in\posop$, then $r\notin\spec{A}$ if and only if
there exist an $\varepsilon>0$ such that
$
 \aleph_{[r-\varepsilon, \,+\infty)}(A)=
\aleph_{[r+\varepsilon, \,+\infty)}(A)
$.
\EOE
\end{rem}

\medskip
\begin{pro}\label{inclusion de espectros}
Let $A\in\posop$. Then, if $\sig{\ese}{A}$ is considered as
acting on $\ese$, it holds
\[\spec{\sig{\ese}{A}}\subseteq \spec{A}.\]
\end{pro}
\bdem
By Proposition \ref{cotita},  if $0 \notin \spec{A}$ then
$0 \notin \spec{\sig{\ese}{A}}$. On the other hand, if $r>0$ and
$r\notin \spec{A}$, then, by Remark \ref {keeping}, there exists
$\varepsilon>0$ such that
$
 \aleph_{[r-\varepsilon, \,+\infty)}(A)=
\aleph_{[r+\varepsilon, \,+\infty)}(A) .
$
Hence,
\[
\aleph_{[r-\varepsilon, \,+\infty)}(\sig{\ese}{A})=P_\ese\wedge
\aleph_{[r-\varepsilon, \,+\infty)}(A)=P_\ese \wedge
\aleph_{[r+\varepsilon, \,+\infty)}(A)= \aleph_{[r+\varepsilon,
\,+\infty)}(\sig{\ese}{A}) .
\]
Thus, $r\notin \spec{\sig{\ese}{A}}$.\edem

\begin{pro}\label{saca nondecreasing}
Let $A\in\posop$, $\ese$ a closed subspace and
$f:[0,+\infty)\rightarrow[0,+\infty)$ a non-decreasing right
continuous function. Then,
\begin{equation}\label{eq saca monotonas dos}
  f(\sig{\ese}{A})=\sig{\ese}{f(A)}
\end{equation}
\end{pro}

\bdem

Given $\la\geq 0$, since $f$ is non-decreasing and right
continuous there exist $\eta\geq 0$ such that
\[
\{\mu:\,f(\mu)\geq \la\}=[\eta,+\infty)
\quad \hbox{ and, for every  }\ \  C \in \cam ,  \quad 
\aleph_{[\la , \infty)}(f(C)) =
\aleph_{[\eta , \infty)}(C).
\]
If $\eta=0$, then $\aleph_{[\la , \infty)}(f(\sig{\ese}{A})) =
\aleph_{[\la , \infty)}(\sig{\ese}{f(A)})=I$. On the other hand,
if $\eta>0$,
\medskip
\begin{align*}
\aleph_{[\la , \infty)}(f(\sig{\ese}{A})) &= \aleph_{[\eta ,
\infty)}(\sig{\ese}{A}) = \aleph_{[\eta, \infty)}(A) \wedge P_\ese
\\&\\&= \aleph_{[\la , \infty)}(f(A)) \wedge P_\ese =\aleph_{[\la ,
\infty)}(\sig{\ese}{f(A)}),
\end{align*}

\smallskip
\noi which shows that $f(\sig{\ese}{A})$ and $\sig{\ese}{f(A)}$
have the same (left) spectral resolution. Hence
$f(\sig{\ese}{A})=\sig{\ese}{f(A)}$\edem

\subsection*{Computation of $\min \spec{\sig{\ese}{A}}.$}

\begin{pro}\label {con lamda}
Let $A \in \posop$. Then, if $\sig{\ese}{A}$ is considered as
acting on $\ese$, it holds
\begin{equation}\label{con la}
 \min \spec{\Csha}
 = \max \{ \la  \ge 0 : A^m \ge  \la ^m P_\ese  , \;
 \forall \ m \in \N \} .
\end{equation}
\end{pro}
\bdem Note that $A^m \ge  \la ^m P_\ese$, $  m \in \N$,
if and only if $\la P_\ese \leqr A$. On the other
hand, since $P_\ese$ and $\sig{\ese}{A} $ commute,
$\la P_\ese \le \sig{\ese}{A} $ \sii
$\la P_\ese \leqr \sig{\ese}{A} $ if and only if $\la P_\ese \in \eme_\rho(\ese,A)$
if and only if $ \la P_\ese \leqr A $.

\edem

\bigskip

\begin{teo}\label{el minimo del espectro}
Let $A\in \posop$. Then, if $\sig{\ese}{A}$ is considered as acting on $\ese$,
\begin{equation}\label{formula para el minimo del espectro}
\barr{rl} \min\spec{\sig{\ese}{A}}&
= \max \{ \la \ge 0: P_\ese \le \aleph _{[\la , \infty )}(A) \} \\&\\
& = \min \{\mu\in\spec{A}:\;R( \aleph _{[\mu , \mu
+ \eps ) } (A) ) \not \inc \cS^\bot\;\;\forall \; \varepsilon>0\}\\&\\
& = \min \{\mu\in\spec{A}:\; P_\ese \  \aleph _{[\mu , \mu
+ \eps ) } (A)  \not = 0\;\;\forall \; \varepsilon>0\}.
\earr
\end{equation}
\end{teo}
\bdem
For any $B \in L(\cS )^+$, $\min \spec{B}
= \max \{ \la \ge 0 :
\aleph _{[\la , \infty)}(B) =I_\ese  \}.
$
Applying this identity to our problem, we get
$\la_0 = \min\spec{\sig{\ese}{A}}= \max \{ \la \ge 0: P_\ese \le \aleph _{[\la , \infty )}(A) \}
$. Then $P_\ese \le \aleph _{[\la_0 , \infty )}(A)$ but $P_\ese \not\le \aleph _{[\la_0+ \eps  , \infty )}(A)$
for every $\eps >0$. So that, $\la_0 \in \{\mu\in\spec{A}:\;P_\ese \  \aleph _{[\mu , \mu
+ \eps ) } (A)  \not = 0 \;\;\forall \; \varepsilon>0\}$,
since if $P_\ese \  \aleph _{[\la_0 , \la_0 + \eps ) } (A)  = 0$, then
$$P_\ese \  \aleph _{[\la_0 + \eps , \infty ) } (A) =P_\ese \ \Big(
\aleph _{[\la_0 , \infty ) } (A) - \aleph _{[\la_0 , \la_0 + \eps ) }(A)\Big)  =
P_\ese \ \aleph _{[\la_0 , \infty ) } (A) = P_\ese ,
$$
i.e. $P_\ese \le \aleph _{[\la_0 + \eps , \infty )}(A)$.
If $\la_0 = 0$, then equation (\ref{formula para el minimo del espectro}) is clear, 
since $[\la_0, \la_0 +\eps)$ is an open subset of $\spec {\sig{\ese}{A}}$.
If $\la_0 > 0$, let  $0\le \la <\la_0$ and $0< \eps < \la_0 -\la$. Then $\la + \eps \le \la_0$. Since
$\la_0 = \max \{ \la \ge 0: P_\ese \le \aleph _{[\la , \infty )}(A) \} $, it holds that
$ P_\ese \aleph _{[\la , \infty ) } (A)
= P_\ese \aleph _{[\la +\eps  , \infty ) } (A) = P_\ese$.
Hence
$$
P_\ese =  P_\ese \aleph _{[\la , \infty ) } (A)  =
 P_\ese \aleph _{[\la , \la+\eps  ) } (A) +P_\ese \aleph _{[\la +\eps, \infty ) } (A)
 = P_\ese \aleph _{[\la , \la+\eps  ) } (A) +P_\ese .
 $$
 Therefore
 $P_\ese \aleph _{[\la , \la+\eps  ) } (A) = 0$, showing equation
 (\ref{formula para el minimo del espectro}).
\edem

\section{The case $\dim \ese=1$.}

\begin{fed}\rm
Suppose that $\dim \ese = 1$ and let $\xi \in \ese$ an unit
vector. For every $A\ge 0$ there exist $\la, \mu  \ge 0$ such that
$\Csha = \la P_\ese$ and  $\csha = \mu P_\ese$. Denote
$\sig{A}{\xi}=\la$ and  $\Sig{A}{\xi}=\mu$.
\end{fed}
\begin{rem}\rm \label{cosas}
Let $\ese $  be the subspace generated by the unit vector $\xi
\in\hil$. There are several ways to compute $\sig{A}{\xi}$ in
terms of $\sig{\ese}{A}$, and similarly $\Sig{A}{\xi}$ in terms of
$\Sig{\ese}{A}$. For example: \ben \item By Theorem \ref{el minimo
del espectro},
\begin{equation} \label{kolmo}
\barr{rl}
\sig{A}{\xi}& = \min \spec{\sig{\ese}{A}} = \min \Big\{ \mu\in\spec{A}: \;
P_\ese  \aleph _{[\mu , \mu + \eps ) } (A)  \not = 0 \; \; \forall \; \varepsilon>0 \Big\}
\\&\\
& = \min \Big\{ \mu\in\spec{A}: \;
\aleph _{[\mu , \mu + \eps ) } (A)\xi  \not = 0 \; \; \forall \; \varepsilon>0 \Big\}
\earr .
\end{equation}
\item By Proposition \ref{con lamda}
$$
\sig{A}{\xi} = \max \{ \la \ge 0 : \api
A^n \eta , \eta \cpi \ge  \la^n |\api \xi , \eta \cpi |^2
\ ,\; \; \forall \; n \in \N , \ \eta \in \hil   \}.
$$
\item
Also $\sig{A}{\xi} = \|\sig{\ese}{A} \xi\| = \api \sig{\ese}{A} \xi , \xi\cpi$. Similar
formulae hold for $\Sig{A}{\xi}$.
\item By Proposition \ref{inclusion de espectros},
$\sig{A}{\xi} \in \spec{A}$. Moreover, by Theorem \ref{trunquetti}
(or Theorem \ref{el minimo del espectro}),
\beq\label{formulae}
\sig{A}{\xi} = \max \{ \la \in \spec{A} : \xi \in R(\aleph_{[\la , \infty )} (A)) \}.
\end{equation}
\een
\end{rem}

\noi
The following result relates the spectral short of operator to one
dimensional subspaces and the spectral order.

\begin{pro} \label{dim1}
Let $A$, $ B \in\cam$. Then $A\leqr B$ \sii $\sig{A}{\xi}\le
\sig{B}{\xi}$ for every unit vector $\xi\in\hil$.
\end{pro}
\bdem One implication follows from Corollary \ref{monotonia}. On
the other hand, suppose that $\sig{A}{\xi}$ $\le \sig{B}{\xi}$ for
every unit vector $\xi\in\hil$. Given $\la\geq 0$ such that 
$\aleph_{[\la , \infty )}({A})\neq 0$, let $\zeta \in
R(\aleph_{[\la , \infty )}({A}))$. By equation (\ref{formulae}),
$\la\le\sig{A}{\zeta}$. Since $\sig{A}{\zeta}\le \sig{B}{\zeta}$,
by equation (\ref{formulae}) we have that $\zeta \in
R(\aleph_{[\la , \infty )}({B}))$. 
Hence $R(\aleph_{[\la , \infty )}({A})) \inc R(\aleph_{[\la , \infty )}({B}))$
for every $\la \ge 0$. By Theorem
\ref{olson},  we deduce that $A\leqr B$. \edem

\begin{pro}\label{robusto}
Let $A \in \posop$ and let $\ese $  be the subspace of $\hil$
generated by the unit vector $\xi$. If $A$ is invertible, then
for $ m \in \N$,
\beq\label{invv}
\Sig{A^{2m}}{\xi}^{1/2m} = \| A ^{-m}\xi\|^{-1/m} = \api A^{-2m}\xi , \xi \cpi ^{-1/2m},
\end{equation}
and
\beq\label{alam}
\sig{A}{\xi} = \lim _{m \to \infty } \|A^{-m} \xi \|
^{-1/m} = \inf _{m \in \Bbb N}  \|A^{-m} \xi \| ^{-1/m}
\end{equation}
If $R(A)$ is closed,  then

\ben
\item  If $\xi \not\in R(A)$, then $\sig{A}{\xi} = 0$.
\item If $\xi \in R(A)$ and $B = A^\dagger$, then
$\sig{A}{\xi} = \lim _{m \to \infty } \|B^m \xi \| ^{-1/m} =
\inf _{m \in \Bbb N}  \|B^{m} \xi \| ^{-1/m}$.
\een

\end{pro}
\bdem

Using Theorem \ref{el minimo del espectro}, the closed range case
easily reduces to the invertible case, by considering $A$ as
acting on  $R(A)$, because $A^\dagger$ acts on $R(A)$ as the
inverse of $A$. Note that, if $R(A)$ is closed, then there exists
$\eps>0$ such that $ \aleph _{[0 ,  \eps ) } (A) = P_{ N  (A)}$.
Therefore $\xi \not\in R(A)$ implies that $P_\ese \aleph _{[0 ,
\eps ) } (A) \not = 0$, and, by Remark \ref{cosas}, that
$\sig{A}{\xi} = 0$.

Suppose that $A$ is invertible. For $m \in \N$, denote by $\eta_m =
A^{-m/2}\xi$. Fix $m \in \N$. By Theorem \ref{propiedades del shorted},
if $\eme _m = A^{-m/2} (\ese ) $, then $\Sig{\ese}{A^m}= A^{m/2} P_{\eme_m} A^{m/2}$,
and
$$\Sig{A^m}{\xi}=\|\Sig{\ese}{A^m} \xi\| = \|A^{m/2} P_{\eme_m} A^{m/2} \xi\|.
$$
Note that $\eme _m $ is  the subspace generated by $\eta_m$, so that $P_{\eme_m} \rho
= \|\eta_m\|^{-2} \api \rho , \eta_m \cpi \eta_m $, $\rho \in \cH$. Then
$$
\barr{rl}
\Sig{A^m}{\xi} = \|A^{m/2} P_{\eme_m} A^{m/2} \ \xi\| & =
\Big\|A^{m/2} \Big( \|\eta_m\|^{-2} \api A^{m/2} \ \xi , \eta_m \cpi \eta_m \Big)\Big\| \\&\\ & =
\|\eta_m\|^{-2} \|\api  \xi , \xi \cpi \ \xi\| = \|\eta_m\|^{-2} . \earr
$$
Therefore $\Sig{A^{2m}}{\xi} = \| A ^{-m} \xi\|^{-2} $, so that
$$\Sig{A^{2m}}{\xi}^{1/2m} = \| A ^{-m}\xi\|^{-1/m} \ , \quad m \in \N .
$$
Equation (\ref{alam}) follows using Remark \ref{cosas} and the definition of $\sig{\ese}{A}$.
\edem

\begin{rem} \label{schurinv}  \rm
Equation (\ref{invv}) and, consequently, Proposition \ref{robusto}, can also be deduced
from the following formula: for every invertible $B \in \cam $ and $\xi \in \hil $ with $\|\xi \| = 1$,
$$\Sig{B}{\xi} = \api B\inv \xi , \xi \cpi \inv .
$$
This formula is the one dimensional case of  the characterization
of Schur complements in terms of the block representation of the
inverse of an operator (see \cite{[pedro]} Lemma 4.7 or, for a
matrix version, Horn-Johnson book \cite{[HoJo]}). \EOE
\end{rem}

\medskip
\noi
Let $A\in \cam$. We shall denote by
$$\specr{A} = \big\{ \sig{A}{\xi} : \|\xi \| = 1\big\} .
$$
By Proposition \ref{inclusion de espectros}, we have that
$\specr {A}\inc \spec{A}$. If $\dim \H < \infty$, it was shown in \cite{[ACS]} 
(see also \cite{[FF]}) that $\specr {A}= \spec{A}$. 
We shall see that this property fails in general. First we fix some notations:
\ben 
\item For $B\in \cam $ we denote by 
$$
\barr{rl}
\specm{A}& = \big\{ \la \in
\spec{A}: \exists \ (\mu_n)_{n\in\N} \ in \ \spec{A} \  such \  that  \
\mu_n > \la \  and \
\mu _n \convsotdpre_{n\rightarrow \infty} \la \big\} \\&\\
& = \big\{ \la \in
\spec{A}: \forall \ \eps >0 \ , \ \aleph_{(\la, \la+ \eps )}(A) \neq 0
\} ,
\earr
$$
i.e. those points  $\la \in \spec{A}$ which are limit point of
$\spec{A} \setminus \{\la\}$ from the right.
\item $\sigma_{pt}(A) = \big\{\la \in \spec{A} : N(A-\la I) \neq \{0\}\big\} $, 
the point spectrum of $A$.
\een

\begin {pro} \label{por izq} Let $A\in \cam$. Then
$$\specr{A} = \specm{A} \cup \sigma_{pt}(A)
=\big\{ \la \in
\spec{A}: \forall \ \eps >0 \ , \ \aleph_{[\la, \la+ \eps )}(A) \neq 0
\}.
$$
In particular, this shows that $\specr{A}$ is dense in
$\spec{A}$.
\end{pro}
\bdem
Let $\la \in \spec{A}$ and let $(\mu_n)_{n\in\N}$ be a sequence  in $\spec{A}$   such   that
$\mu _n \convsotdpre_{n\rightarrow \infty} \la$. Denote by $\la _0 = \mu_1+1$ and
$\la_n = \frac{1}{2}(\mu _{n+1}+\mu_n)$,
$n \in \N$. Note that, since $\mu_n \in (\la_n, \la_{n-1})$,
then $\aleph_{(\la_n, \la_{n-1})}(A) \neq 0$.
We take, for every $n \in \N$, an unit vector
$\xi_n \in R(\aleph_{(\la_n, \la_{n-1})}(A)).$
Consider the unit vector
$$
\xi = \sum_{n \in \N} \ \frac{\xi_n}{2^n}  \ .
$$
Recall formula (\ref{formulae}), wich says that
$\sig{A}{\xi} = \max \{ \mu \in \spec{A} : \xi \in R(\aleph_{[\mu , \infty )} (A)) \}.$
It is clear by construction of $\xi$  that in our case we get $\sig{A}{\xi} =\la$, because
$\la = \inf_n \mu_n = \inf_n \la_n $.
If $\la \in \sigma_{pt}(A)$, just take $\xi \in  N   (A-\la I)$ and clearly $\sig{A}{\xi}=
\Sig{A}{\xi}=\la$.

Now suppose that $\la \in \spec{A}$ but $\la \notin \specm{A} \cup \sigma_{pt}(A)$. This means that
there exists $\eps>0$ such that $\aleph_{[\la, \la+ \eps )}(A) = 0$. Therefore,
for any unit vector $\xi$, it is impossible that
$$
\la = \max \{ \mu \in \spec{A} : \xi \in R(\aleph_{[\mu , \infty )} (A)) \},
$$
because if $\xi \in R(\aleph_{[\la , \infty )} (A))$, then
$\xi \in R(\aleph_{[\la +\eps  , \infty )} (A))$.
\edem
\medskip
\begin{rem}\rm
If $A\in \cam $ is not invertible, then $0 \in \spec {A}$. 
If $0$ were an isolated point of $\spec {A}$ then $A$ must 
have closed range. So that, $ N (  A) \neq \{0\}$.
Otherwise $\aleph_{(0, \eps)} (A) \neq 0 $ for every $\eps >0$. 
This shows that $0 \in \specr{A}$. 
More generally, for $A \in \cam$, it holds that $\la_{min}(A) = \min \spec{A} \in \specr{A}$.
On the other hand, by Proposition \ref{por izq}, $\|A\|\in \specr{A} $ \sii 
$\|A\|$ is an eigenvalue of $A$.
\end{rem}\
\medskip
\begin{rem}\rm
For $A \in \cam$, we shall denote by $R_0(A)$ the subspace
$$
R_0(A) =\bigcup_{\la > 0} R(\aleph_{[\la , \infty)} (A)) .
$$
If $R(A)$ is closed, then $R_0(A) = R(A)$, since $0$ is an isolated point of $\spec{A}$.
But in other case, $R_0(A)$ is properly included in $R(A)$, but it is
still a dense subspace of $\overline{R(A)}$.
We are interested in this subspace
because, by formula (\ref{formulae}), if $\xi \in \H$ an unit vector, then
$\sig{A}{\xi} \neq 0$ \sii $\xi \in R_0(A)$. 
\EOE
\end{rem}

\subsection{Kolmogorov's complexity}

Given an invertible matrix  $A\in L(\C^m )^+ $ and $\xi \in \C^m$ a unit vector,
J. I. Fujii and M. Fujii
\cite{[FF]} define the Kolmogorov's complexity:
\beq\label{kol1}
\kol{A}{\xi} = \lim _{n \to \infty} \frac{\log (\pint {A^n \xi , \xi})}{n} =
\log \lim _{n \to \infty} \pint {A^n \xi , \xi}  ^{1/n} .
\end{equation}
Using formula (\ref{invv}), we can see that the limit is, in fact, a supremum; and we have
the identity
\beq\label{kol2}
\kol{A}{\xi} = \log \sig{A^{-1/2}}{\xi}^{-2} = \log \sig{A^{-1}}{\xi } \inv .
\end{equation}
This shows, using formulae (\ref{kolmo}) and (\ref{formulae}), the following formula:
\beq\label{kol3}
\barr{rl}
\exp \kol{A}{\xi} &= \min \big\{ \la \in \spec{A} : \xi \in R(\aleph_{(-\infty , \la]} (A))\big\}
\\&\\
& = \max \Big\{ \mu\in\spec{A}: \;
\aleph _{( \mu - \eps, \mu] } (A)\xi  \not = 0 \; \; \forall \; \varepsilon>0 \Big\} \  .
\earr
\end{equation}
With these identities in mind we generalize the notion of
Kolmogorov's complexity in two directions: firstly we define it
for infinite dimensional Hilbert spaces; secondly we remove the
hypothesis of invertibility of $A$. Note that the own notion of
spectral shorted operator is, in some sense, a generalization of
the Kolmogorov's complexity relative to arbitrary (not necesarily
onedimensional) closed subspaces of a Hilbert space $\H$.

If $\H$ is a Hilbert space and $A\in \cam$ is invertible, then we
just have to define $\kol{A}{\xi} $ as in equation (\ref{kol2})
or, equivalently, equation (\ref{kol3}). It is easy to see that
this is equivalent to define it as in the finite dimensional
setting, as in equation (\ref{kol1}). We should mention that some
of the properties of $\kol {A} {\xi}$ proved by J. I. Fujii and M.
Fujii fail if $\hil$ is infinite dimensional. As an example, the
identity
$$
\sigma (A) = \big\{ \exp (\kol {A} {\xi}) : \|\xi \| = 1 \big\} .
$$

\begin{fed}\rm
Given $\xi \in \H$ and $A\in \cam$,
we denote by
$$
\kkol{A}{\xi} = \lim _{n \to \infty} \pint {A^n \xi , \xi}  ^{1/n} ,
$$
that is, $\kkol{A}{\xi} = \exp \kol{A}{\xi}$ in the cases where
$\kol{A}{\xi}$ is defined. 
\end{fed}
\begin{rem}\label{kolchico}\rm
If $\xi \in \H$ and $A\in \cam$, then
\ben
\item If $\|\xi\| = 1$, then the sequence $\pint {A^n \xi , \xi}  ^{1/n}$ is increasing.
So that, $\lim _{n \to \infty} \pint {A^n \xi , \xi}  ^{1/n}$  exists
for every $\xi \in \H$.
\item $\kkol{A}{\xi} = \kkol{A}{ a \xi}$ for every $0 \neq a \in \C$.
\item $\kkol{A}{\xi} = \kkol{A}{\aleph_{[\la , \infty)} (A) \xi}$ for every
$\la >0 $ such that $\aleph_{[\la , \infty)} (A)\xi \neq 0$.
\een
Indeed,  by H\"older inequality for states (also by Jensen inequality, see \cite{[AMS]}),
if $\|\xi \| = 1$, $p \ge 1$ and $1/p + 1/q =  1$, then
$$\pint {A^p \xi , \xi}  ^{1/p} \pint {I^q \xi , \xi}  ^{1/q} = \pint {A^p \xi , \xi}  ^{1/p}
\ge \pint{A\xi , \xi}.
$$
Applying this inequality to
$A^{n}$ with $p = (n+1)/n$ one gets that $\pint {A^n \xi , \xi}  ^{1/n} \le
\pint {A^{n+1} \xi , \xi}  ^{1/n+1}$. \\
Item 2 follows from the fact that $|a|^{2/n} \conv 1$. To show  3, suppose that
$\|\xi \| = 1$ and denote by
$\xi_1 = \aleph_{[\la , \infty)} (A)\xi $ and $\xi_2 = \xi -\xi_1$. Then,
since $\aleph_{[\la , \infty)} (A)$ commutes with $A$, for every $n \in \N$,
$$
\barr{rl}
\pint {A^n \xi_1 , \xi_1} & \le \pint {A^n \xi_1 , \xi_1}  +
\pint {A^n \xi_2 , \xi_2} =
\pint {A^n \xi , \xi}
\\&\\
&
\le \pint {A^n \xi_1 , \xi_1} +\la^n
 \le
(1+ \|\xi_1\|^{-2} ) \pint {A^n \xi_1 , \xi_1} .  \earr
$$
This shows that  $\kkol{A}{\xi} = \kkol{A}{\xi_1 }$, since
$(1+ \|\xi_1\|^{-2} )\ ^{1/n} \conv 1$.
\EOE
\end{rem}

\noi
Recall that, for $A \in \cam$, we denote by
$
R_0(A) =\bigcup_{\la > 0} R(\aleph_{[\la , \infty)} (A)) .
$

\begin{pro}\label{denso} Let $A\in \cam $ and $0 \neq \xi \in \H$.
Then $\kkol{A}{\xi} \neq 0$ \sii $P_{\overline{R(A)}}\  \xi \in R_0(A)$. Moreover,
equation (\ref{kol3}) holds in general:
\beq\label{kol n}
\barr{rl}
\kkol{A}{\xi} &= \min \big\{ \la \in \spec{A} : \xi \in R(\aleph_{(-\infty , \la]} (A))\big\}
\\&\\
& = \max \Big\{ \mu\in\spec{A}: \;
\aleph _{( \mu - \eps, \mu] } (A)\xi  \not = 0 \; \; \forall \; \varepsilon>0 \Big\} 
\\&\\
& = \sup \Big\{ \mu\in\spec{A}: \; \aleph_{[\mu , \infty)} (A)\xi \neq 0 \Big\} \ .
\earr
\end{equation}
\end{pro}
\bdem
Let $\la = \sup \Big\{ \mu\in\spec{A}: \; \aleph_{[\mu , \infty)} (A)\xi \neq 0 \Big\}$.\\
If $\mu > \la$, then $\xi \in R (\aleph_{(- \infty , \mu ]} (A))$, so that
$ \pint {A^n \xi , \xi} \le \mu^n \|\xi\|^2 $ for $n \in \N$, and $\kkol{A}{\xi} \le \mu$.
On the other hand, if $\mu < \la$ then  $\aleph_{[\mu , \infty)} (A)\xi = \xi_1 \neq 0 $,
and, by Remark \ref{kolchico},
$\kkol{A}{\xi} = \kkol{A}{\xi_1 } \ge \mu$, since $\pint {A^n \xi_1 , \xi_1}
\ge \mu^n \|\xi_1\|^2 $ for every $n \in \N$.
This shows that $\kkol{A}{\xi} = \la$. The other equalities are straightforward, 
by spectral theory.
\edem

\noi
By Proposition \ref{inclusion de espectros}, we have that
$\specr {A}\inc \spec{A}$ and, therefore, if $A$ is invertible,
$$
\big\{  \kkol {A}{\xi} : \|\xi \|\neq 0 \big\} =
\big\{ \sig{A\inv}{\xi}\inv : \|\xi \| = 1 \big\} \inc \spec{A\inv}\inv = \spec{A}.
$$
As we shall see below, the reverse inclusion fails in general:

\begin{pro}\label{por der}
Let $A>0$. Then
$$
\barr{rl}
\big\{  \kkol {A}{\xi} : \|\xi \|\neq 0 \big\} & =
\ \specM{A} \cup \sigma_{pt}(A) \\&\\
& = \ \big\{ \la \in
\spec{A}:    \aleph_{(\la+ \eps, \la ]}(A) \neq 0 \ , \ \forall \ \eps >0 \
\},
\earr
$$
where $ \specM{A}$ is the set of points in $\spec{A}$ which are limit point of
$\spec{A} \setminus \{\la\}$ from the
left. The set $\big\{  \kkol {A}{\xi} : \|\xi \|= 1 \big\}$
is also dense in $\spec{A}$.
\end{pro}
\bdem
It is a consequence of Proposition \ref{por izq} (applied to $A\inv$) and the identity
$$
\big\{  \kkol {A}{\xi} : \|\xi \| \neq 0 \big\} =
\big\{  \kkol {A}{\xi} : \|\xi \|= 1 \big\} =
\big\{ \sig{A\inv}{\xi}\inv : \|\xi \| = 1 \big\} .
$$
\edem
\begin{rems}\rm \
\ben
\item Proposition \ref{por der} is also true for a general $A \in \cam$.
The proof is similar to the proof of Proposition \ref{por izq},
by using equation (\ref{kol n}) instead of (\ref{formulae}).
\item
Let $\H = \ell^2 (\N )$, denote by $e_n $, $n \in \N$, the cannonical orthonormal basis of $\H$,
and consider the diagonal invertible operators $A, B \in \cam$ defined by
$$
A(e_n ) = \big(2 + \frac{1}{n}\big) e_n  \quad , \quad
B(e_n ) = \big(2 - \frac{1}{n}\big) e_n \ , \quad n \in \N .
$$
It is easy to see, using Propositions \ref{por izq} and \ref{por der}, that
$2 \notin \big\{ \kkol {A}{\xi} : \|\xi \|= 1 \big\}$ and $2 \notin \specr{B}$.
\item
If $C\in \cam$, then  $\|C\| \in \big\{ \kkol {C}{\xi} :
\|\xi \|= 1 \big\}$ and $\la_{min}(C) \in \specr{C}$. On the other hand,
if $A$ and $B$ are the operators of the previous example,
$\| B\| = 2 \notin \specr {B}$ and $\la_{min}(A) = 2  \notin \big\{ \kkol {A}{\xi} :
\|\xi \|= 1 \big\}$.
\een
\end{rems}

\begin{rem}[Operators with closed range] \rm
Suppose that $A\in \cam$ and  $R(A)$ is closed.
\ben
\item
 If $\xi \in R(A)$ is an unit vector, then, by Proposition
\ref{robusto}, $\kkol{A}{\xi} = \sig{A^\dag}{\xi}\inv$.
\item If $\xi \notin  N  ( A)$ and $\xi \notin R(A)$, then
the behaviors of the Kolmogorov complexity and the spectral shorted operator, relative to $\xi$
are different. By Proposition \ref{robusto}, $\sig{A}{\xi} = \sig{A^\dag}{\xi} = 0$.
On the other hand, if $P = P_{R(A)}$, then $P\xi \neq 0$ and
$$
\kkol{A}{\xi} = 
\lim _{n \to \infty} \pint {A^n P \xi , P \xi}  ^{1/n} =
 \kkol{A}{\frac{P\xi}{\|P\xi\|}} =  \sig {A^\dag }{\frac{P\xi}{\|P\xi\|}}\inv  \neq 0 .
 $$
\een
\end{rem}


\begin{thebibliography}{XXXXXX}
\fontsize {10}{10}
\bibitem{[AT1]}W. N. Anderson, Shorted operators, SIAM J.
Appl. Math. 20 (1971), 520-525.
\bibitem{[AT2]}W. N. Anderson and  G. E. Trapp,
Shorted operators II, SIAM J. Appl. Math. 28 (1975), 60-71.
\bibitem{[ACS]} J. Antezana, G. Corach and D. Stojanoff, Spectral shorted matrices,
Linear Algebra Appl. {381} (2004), 197-217.
\bibitem{[AMS]} J. Antezana, P. Massey and D. Stojanoff, Jensen inequalities and
majorization, preprint.
\bibitem{[FF]} Jun Ichi Fujii and  Masatoshi Fujii,
Kolmogorov's complexity for positive definite matrices, Special
issue dedicated to Professor T. Ando. Linear Algebra Appl. 341
(2002), 171-180.
\bibitem{[HoJo]} R.A. Horn and C.R. Johnson, Matrix
analysis, Cambridge University Press, 1985.
\bibitem{[KR]} R. V. Kadison and J. R. Ringrose,
Fundamentals of the theory of operator algebras. Vol. I.
Pure and Applied Mathematics, 100. Academic Press, Inc., New York, 1983.
\bibitem{[K]} M. G. Krein,  The theory of self-adjoint
extensions of semibounded Hermitian operators and its
applications, Mat. Sb. (N. S.) 20 (62) (1947), 431-495
\bibitem{[pedro]} P. Massey and D. Stojanoff, Generalized Schur Complements
and P-com\-plementable Operators, Linear Alg. Appl. (to appear).
\bibitem{[mathias]} Chi-Kwong Li and R. Mathias, Extremal characterizations
of the Schur complement and resulting inequalities, SIAM Rev. 42
(2000), no. 2, 233--246.
\bibitem{[Ols]} M. P. Olson, The selfadjoint operators of
a von Neumann algebra form a conditionally complete lattice, Proc.
Amer. Math. Soc., 28 (1971) 537-544.
\bibitem{[Ped]} G. K. Pedersen, $C^{*}$-algebras and their
automorphism groups. London Mathematical Society Monographs, 14.
Academic Press, London-New York, 1979.
\bibitem{[Ped2]} G. K. Pedersen, Analysis now.
Graduate Texts in Mathematics, 118. Springer-Verlag, New York, 1989.
\bibitem{[PEKA]}E. L. Pekarev, Shorts of operators and some
extremal problems, Acta Sci. Math. (Szeged) 56 (1992), 147-163.
\bibitem{[PS]} E. L. Pekarev and J. L. Smul'jan, Parallel
addition and parallel substraction of operators, (Russian,
English) Math. USSR, Izv. 10(1976), 351-370
\bibitem{[smu]} J. L. Smul'jan,  A Hellinger operator integral.
(Russian) Mat. Sb. (N.S.) 49 (91) 1959 381-430. English transl.
AMS Transl. 22 (1962), 289-337.
\end{thebibliography}
\end{document}